\documentclass[12pt,reqno]{amsart}
\usepackage{amsmath,amstext,amssymb,amscd,hyperref}
\usepackage{verbatim}
\usepackage{enumerate}
\usepackage{mathrsfs}
\usepackage[dvipsnames,usenames]{xcolor}

\usepackage{amsthm}

\newtheorem{theorem}{Theorem}[section]
\newtheorem{lemma}[theorem]{Lemma}

\theoremstyle{definition}
\newtheorem{definition}[theorem]{Definition}

\theoremstyle{remark}

\numberwithin{equation}{section}
\author[C. Cao]{Chongsheng Cao}
\address{Department of Mathematics \& Statistics \\Florida International University\\
Miami, Florida 33199, USA} \email{caoc@fiu.edu}

\author[Y. Guo]{Yanqiu Guo}
\address{Department of Mathematics \& Statistics \\Florida International University\\
Miami, Florida 33199, USA} \email{yanguo@fiu.edu}

\author[E. S. Titi]{Edriss S. Titi}
\address{Department of Mathematics \\ Texas A\&M University\\
College Station, TX 77843, USA \textbf{AND} Department of Computer Science and Applied Mathematics \\ Weizmann Institute of Science \\ Rehovot 7610001 Israel}
 \email{titi@math.tamu.edu, edriss.titi@weizmann.ac.il}

\title[The Hasegawa-Mima model with partial dissipation]
{Global strong solutions for the three-dimensional Hasegawa-Mima model with partial dissipation}
\date{January 6, 2018}

\begin{document}

\maketitle

\begin{abstract}
We study the three-dimensional Hasegawa-Mima model of turbulent magnetized plasma with horizontal viscous terms and a weak vertical dissipative term. 
In particular, we establish the global existence and uniqueness of strong solutions for this model.
\end{abstract}

\section {Introduction} \label{S1}

\subsection{Literature}   \label{literature}
In 1977, Hasegawa and Mima introduced a system in \cite{HM, HM2} to elucidate the drift wave turbulence in Tokamak, the most advanced magnetic confinement device. The three-dimensional inviscid Hasegawa-Mima equations can be written as (cf. \cite{ANP, CFT, HM, HM2, MH, ZG}):
\begin{align}  
&\frac{\partial w}{\partial t} + J(\phi, w) + \frac{\partial \phi}{\partial z}=0,  \label{HM-1} \\
&\frac{\partial}{\partial t} (\Delta_h \phi - \phi) + J(\phi, \Delta_h \phi) + \gamma \frac{\partial \phi}{\partial y} - \frac{\partial w}{\partial z} = 0,  \label{HM-2}
\end{align}
where $J(f,g) = \frac{\partial f}{\partial x}   \frac{\partial g}{\partial y} -  \frac{\partial f}{\partial y} \frac{\partial g}{\partial x} $ is the Jacobian 
and $\Delta_h=\frac{\partial^2}{\partial x^2}+ \frac{\partial^2}{\partial y^2}$ is the horizontal Laplacian. System (\ref{HM-1})-(\ref{HM-2}) describes the coupling of the drift modes to the ion-acoustic waves that propagate along the magnetic field. Here, $\phi$ is the electrostatic potential, and simultaneously is the stream function for the horizontal flow in the $xy$-plane. Moreover, $w$ represents the normalized ion velocity in the $z$-direction, and $\gamma$ is a constant which is proportional to the density gradient. 

Like the three-dimensional Euler equations of inviscid incompressible fluid, the only conserved quantity for the 3D Hasegawa-Mima equations (\ref{HM-1})-(\ref{HM-2}) is the kinetic energy, and the global regularity problem is open. Nevertheless, by adding the full viscosity to (\ref{HM-1})-(\ref{HM-2}), Zhang and Guo \cite{ZG} proved the global regularity and the existence of global attractors for a viscous and forced 3D Hasegawa-Mima model using standard tools from the theory of Navier-Stokes equations. On the other hand, Cao, Farhat and Titi \cite{CFT} proposed and studied an inviscid three-dimensional modified version of (\ref{HM-1})-(\ref{HM-2}), the pseudo-Hasegawa-Mima equations:
\begin{align}
&\frac{\partial w}{\partial t} + \mathbf u \cdot \nabla_h w - U_0 L \frac{\partial \omega}{\partial z}=0,  \label{PHM-1}\\
&\frac{\partial \omega}{\partial t} + \mathbf u \cdot \nabla_h \omega - \frac{U_0}{L} \frac{\partial w}{\partial z}=0,   \label{PHM-2} 
\end{align}
with $\nabla_h \cdot \mathbf u=0$, for some constant $U_0$, where $\textbf u=(u,v)^{tr}$ is the horizontal component of the velocity vector field $(u,v,w)^{tr}$, and 
$\omega=\nabla_h \times \mathbf u$ is the vorticity. The operator $\nabla_h=(\frac{\partial}{\partial x}, \frac{\partial}{\partial y})^{tr}$ is the horizontal gradient. In particular, the global well-posedness of the weak solutions to (\ref{PHM-1})-(\ref{PHM-2}) was established in \cite{CFT}. Observe that $\omega$ in (\ref{PHM-1})-(\ref{PHM-2}) plays the role of the term $\Delta_h \phi-\phi$ in (\ref{HM-1})-(\ref{HM-2}). Therefore, system (\ref{PHM-1})-(\ref{PHM-2}) is a modified version of the Hasegawa-Mima equations (\ref{HM-1})-(\ref{HM-2}), with the essential difference that the term $\frac{\partial \phi}{\partial z}$ is replace by $\frac{\partial \omega}{\partial z}$. Nevertheless, model (\ref{PHM-1})-(\ref{PHM-2}) is simpler than (\ref{HM-1})-(\ref{HM-2}) in the sense that it has a nice mathematical structure. Indeed, adding and subtracting (\ref{PHM-1}) and (\ref{PHM-2}) yield a three-dimensional coupled transport system with collinear transport velocities in opposite directions leading to an intensified shear in the vertical direction, which results in exponential growth in the relevant estimates for (\ref{PHM-1})-(\ref{PHM-2}) in \cite{CFT}. 

It is worth mentioning other interesting models describing plasma turbulence. For instance, Hasegawa and Wakatani proposed equations for a two-fluid model which describe the resistive drift wave turbulence in Tokamak (cf. \cite{HW, HW2}). The existence and uniqueness of strong solutions to the Hasegawa-Wakatani equations have been established by Kondo and Tani \cite{KT}.

In the context of geophysical fluid dynamics, there are certain models resemble the structure of Hasegawa-Mima equations (\ref{HM-1})-(\ref{HM-2}). In particular, Charney \cite{Charney} and Obukhov \cite{Obuhov} derived the following two-dimensional shallow water model from the Euler equations with free surface under a quasi-geostrophic velocity field assumption:
\begin{align}  \label{CO}
\frac{\partial}{\partial t}(\Delta_h \phi_0-F\phi_0)  + J(\phi_0,\Delta_h \phi_0) + J(\phi_0,\phi_B + \beta y)=0.
\end{align}
Here $\phi_0(x,y)$ is the amplitude of the surface perturbation at the lowest order in the Rossby number, and the equation $z=\phi_B(x,y)$ describes the given bottom topography. $F$ is the Froude number. One may refer to \cite{Pedlosky} for a derivation of model (\ref{CO}). For the simple case when $\phi_B$ is a constant representing a flat bottom, (\ref{CO}) reduces to the Hasegawa-Mima-Charney-Obukhov equation:
\begin{align}  \label{2DHM}
\frac{\partial}{\partial t}(\Delta_h \phi_0-F\phi_0)  + J(\phi_0,\Delta_h \phi_0) + \beta \frac{\partial \phi_0}{\partial x}=0.
\end{align}
Since (\ref{2DHM}) bears a close resemblance to the two-dimensional Euler equations, the standard tools for handling the 2D Euler equations can be adopted to analyze (\ref{2DHM}). Indeed, Guo and Han \cite{GH} proved the global existence and uniqueness of solutions for (\ref{2DHM}). For other results concerning (\ref{2DHM}) see, e.g., Paumond \cite{P}, and Gao and Zhu \cite{GZ}.

It is worth mentioning that one may refer to the monographs \cite{Majda,Pedlosky} as well as the papers \cite{JKMW,JK} for other relevant geophysical models.

\vspace{0.1 in}

\subsection{The model}
Motivated by the Hasegawa-Mima equations and the Charney-Obukhov equations mentioned in subsection \ref{literature}, we introduce and study in this paper the following three-dimensional Hasegawa-Mima model with horizontal viscous terms and a weak vertical dissipative term:
\begin{align}    
&\frac{\partial w}{\partial t} + \mathbf u \cdot \nabla_h w - \frac{\partial \psi}{\partial z} = \frac{1}{Re} \Delta_h w,   \label{model-1} \\
&\frac{\partial \omega}{\partial t} + \mathbf u \cdot \nabla_h \omega - \frac{\partial w}{\partial z} = \frac{1}{Re} \Delta_h \omega  + \epsilon^2 \frac{\partial^2\psi}{\partial z^2},   \label{model-2}\\
&\nabla_h \cdot \mathbf u=0.   \label{model-4}
\end{align}  
The velocity vector field $(u,v,w)^{tr}$ defined in $\Omega= [0,L]^2 \times [0,1]$ satisfies the periodic boundary condition with the horizontal velocity $\textbf u=(u,v)^{tr}$. The stream function $\psi$ for the horizontal flow is defined as $\psi= (-\Delta_h)^{-1} \omega$ with $\int_{[0,L]^2} \psi dx dy =0$, and $\omega=\nabla_h \times \mathbf u$. We denote $\nabla_h=\left(\frac{\partial }{\partial x},  \frac{\partial }{\partial y}\right)^{tr}$ and $\Delta_h=\frac{\partial^2}{\partial x^2}+\frac{\partial^2}{\partial y^2}$. The constant $Re$ is the Reynolds number.

System (\ref{model-1})-(\ref{model-4}) bears a resemblance as the three-dimensional Hasegawa-Mima equations (\ref{HM-1})-(\ref{HM-2}) with the difference that Hasegawa-Mima equations are inviscid, whereas model (\ref{model-1})-(\ref{model-4}) is regularized by horizontal viscosity and a partial vertical dissipation. The purpose of introducing and investigating (\ref{model-1})-(\ref{model-4}) is to shed light on the analysis of the inviscid Hasegawa-Mima equations (\ref{HM-1})-(\ref{HM-2}).

Mathematically, the difficulty of establishing the global regularity for system (\ref{model-1})-(\ref{model-4}) lies in the following aspects: 
\begin{enumerate} [(i)]
\item The physical domain is three-dimensional.
\item The regularizing viscosity acts only on the horizontal variables.
\item The system contains the troublesome term $\frac{\partial \psi}{\partial z}$.
\end{enumerate}

Since the lack of the viscosity in the vertical direction provides great challenge for establishing the global regularity, we impose a weak dissipative term $\epsilon^2 \frac{\partial^2 \psi}{\partial z^2}$ in the equation (\ref{model-2}). Since $\psi= (-\Delta_h)^{-1} \omega$, we remark that, as a dissipation, $\frac{\partial^2 \psi}{\partial z^2}$ is weaker than the vertical viscosity $\frac{\partial^2 \omega}{\partial z^2}$. In a priori estimates conducted in section \ref{sec-est}, the dissipative term $\epsilon^2 \frac{\partial^2 \psi}{\partial z^2}$ plays a vital role in controlling the terms $-\frac{\partial w}{\partial z}$ and $-\frac{\partial \psi}{\partial z}$ with the help of an anisotropic Ladyzhenskaya type inequality (see Lemma \ref{lemma1}).

\vspace{0.1 in}

\subsection{Preliminaries}
In this subsection, we introduce some preliminaries that will be used later in our analysis. Recall the three-dimensional periodic space domain $\Omega=[0,L]^2 \times [0,1]$. Throughout, the norm for the $L^p(\Omega)$ space, for $p\in [1,\infty]$, is denoted by $\|f\|_p$. The inner product of $f$ and $g$ in the $L^2(\Omega)$ space is denoted by 
$(f,g)=\int_{\Omega} fg dx dy dz$. As usual, the Sobolev space $H^1(\Omega)=\{f\in L^2(\Omega): \nabla f \in L^2(\Omega)\}$. In addition, we define the following Hilbert space:
\begin{align}   \label{iden-8}
H^1_h(\Omega)= \left\{f\in L^2(\Omega):   \nabla_h f \in L^2(\Omega) \right\},
\end{align}
that features the inner product $(f,g)_{H^1_h(\Omega)}=(f,g)+(\nabla_h f, \nabla_h g)$.

For sufficiently smooth functions $f$, $g$ and $\mathbf u$, with $\nabla_h \cdot \mathbf{u}=0$, integration by parts yields
\begin{align}   \label{iden-0}
\left(\mathbf{u} \cdot \nabla_h f,g\right)=- \left(\mathbf{u} \cdot \nabla_h g,f\right),
\end{align}
which immediately implies that  
\begin{align}   \label{iden-1}
\left(\mathbf{u} \cdot \nabla_h f,f\right)=0.
\end{align}
Recall that the horizontal velocity $\mathbf u$, the vertical vorticity $\omega$, and the stream function $\psi$ for the horizontal flow have the following relations:
\begin{align}   \label{iden-4}
\omega=\nabla_h \times \mathbf u=v_x-u_y, \;\;\; \omega=-\Delta_h \psi,\;\;\;  \mathbf u = (\psi_y,-\psi_x)^{tr},
\end{align}
where $\int_{[0,L]^2} \psi dx dy=0$.
It follows that, if $\omega \in L^2(\Omega)$, then
\begin{align}  \label{iden-3}
(\omega,\psi)=\|\mathbf u\|_2^2.
\end{align}
In addition, for sufficiently smooth functions $f$, $\textbf u$ and $\psi$ such that $\mathbf u = (\psi_y,-\psi_x)^{tr}$, observe that $\mathbf u \cdot \nabla_h  \psi=0$, then apply (\ref{iden-0}) to deduce 
\begin{align}   \label{iden-2}
\left(\mathbf u \cdot \nabla_h  f, \psi\right) = -\left(\mathbf u \cdot \nabla_h  \psi, f\right)=0.
\end{align}

\vspace{0.1 in}

\subsection{Main result}
Before we state the main result of the paper, we give a definition of a strong solution for system (\ref{model-1})-(\ref{model-4}).
\begin{definition}  \label{def-sol}
We call $(\mathbf u,w)^{tr}=(u,v,w)^{tr}$ a \emph{strong solution} on $[0,T]$ for system (\ref{model-1})-(\ref{model-4}) if 
\begin{enumerate}[(i)]
\item  $(\mathbf u,w)^{tr}$ has the following regularity:
\begin{align}   \label{main-2}
\begin{cases}
\mathbf u, \, w \in L^{\infty}(0,T;H^1(\Omega))\cap C([0,T];L^2(\Omega));\\
\Delta_h \mathbf u, \,  \Delta_h w, \, \omega_z,  \, \nabla_h w_z,  \,  \psi_{zz}      \in L^2(\Omega \times (0,T));   \\
\mathbf u_t,  \,   w_t  \in L^2(\Omega \times (0,T));  
\end{cases}
\end{align}
\item the equations below hold in the following sense:
\begin{align*}    
&\frac{\partial w}{\partial t} + \mathbf u \cdot \nabla_h w - \frac{\partial \psi}{\partial z} = \frac{1}{Re} \Delta_h w,     \;\;\text{in}\;\; L^2(\Omega \times (0,T));\\
&\frac{\partial \omega}{\partial t} + \mathbf u \cdot \nabla_h \omega - \frac{\partial w}{\partial z} = \frac{1}{Re} \Delta_h \omega  + \epsilon^2 \frac{\partial^2\psi}{\partial z^2},
\;\;\text{in}\;\; L^2(0,T; H^1_h(\Omega)'),
\end{align*}    
with $\nabla_h \cdot \mathbf u =0$, where $\omega = \nabla_h \times \mathbf u$, $\psi=(-\Delta_h)^{-1}\omega$ with $\int_{[0,L]^2}\psi dx dy =0$, and $(H^1_h(\Omega))'$ is the dual of the space $H^1_h(\Omega)$, defined in (\ref{iden-8}).
\end{enumerate}
\end{definition}

Now we are ready to state the main result of the paper: the global existence, uniqueness, and continuous dependence on initial data of strong solutions for our model (\ref{model-1})-(\ref{model-4}).
\begin{theorem}  \label{main}
Let $T>0$. Assume $(\mathbf u_0, w_0)^{tr} \in (H^1(\Omega))^3$, then system (\ref{model-1})-(\ref{model-4}) admits a unique strong solution $(\mathbf u,w)^{tr}$ on $[0,T]$ in the sense of Definition \ref{def-sol} satisfying the initial condition
$(\mathbf u(0),w(0))^{tr}= (\mathbf u_0, w_0)^{tr}$. 
Moreover, the energy equality is valid for every $t\in [0,T]$:
\begin{align}    \label{energy}
&\frac{1}{2}\left(\|w(t)\|_2^2 + \|\mathbf u(t)\|_2^2 \right) 
+\int_0^t \left[ \frac{1}{Re} \left(\|\nabla_h w\|_2^2 + \|\nabla_h \mathbf u\|_2^2 \right)  + \epsilon^2 \|\psi_z\|_2^2 \right]ds  \notag\\
&=\frac{1}{2}\left(\|w_0\|_2^2 + \|\mathbf u_0\|_2^2  \right).
\end{align}
In addition, the $H^1(\Omega)$ norm of the solution $(\mathbf u,w)^{tr}$ has a uniform bound independent of $T$. That is, 
\begin{align*}
\sup_{0\leq t \leq T} \left(\|\mathbf u(t)\|_{H^1(\Omega)}^2 + \|w(t)\|_{H^1(\Omega)}^2\right) \leq K,
\end{align*}
where $K$ is independent of $T$, but depends only on $Re$, $\epsilon$, $L$, $\|\mathbf u_0\|_{H^1(\Omega)}$ and $\|w_0\|_{H^1(\Omega)}$. 
Furthermore, if $\{(\mathbf u_0^n,w_0^n)^{tr}\}$ is a bounded sequence of initial data in $H^1(\Omega)$ such that $(\mathbf u_0^n,w_0^n)^{tr} \rightarrow (\mathbf u_0,w_0)^{tr}$ in $L^2(\Omega)$, then the corresponding strong solutions $(\mathbf u^n,w^n)^{tr}$ and $(\mathbf u,w)^{tr}$ satisfy
$(\mathbf u^n,w^n)^{tr} \rightarrow (\mathbf u,w)^{tr}  \;\; \text{in} \;\; C([0,T];L^2(\Omega)).$
\end{theorem}

\vspace{0.2 in}

\section{{\it A priori} estimates}   \label{sec-est}
In this section, we assume that system (\ref{model-1})-(\ref{model-4}) holds for smooth functions and we establish the following formal {\it a priori} estimates. However, as we will show in section \ref{sec-existence}, these formal estimates can be justified rigorously by establishing them first for the Galerkin approximation system and then passing to the limit using the appropriate Aubin compactness theorem. 

\subsection{Estimate for $\|w\|_2^2 + \|\mathbf u\|_2^2$} Taking the $L^2(\Omega)$ inner product of the system (\ref{model-1})-(\ref{model-2}) with $(w,\psi)^{tr}$ yields
\begin{align}    \label{L2-1}
\frac{1}{2}\frac{d}{dt} \left(\|w\|_2^2 + \|\mathbf u\|_2^2 \right) +\frac{1}{Re} \left(\|\nabla_h w\|_2^2 + \|\nabla_h \mathbf u\|_2^2 \right) + \epsilon^2 \|\psi_z\|_2^2  = 0,
\end{align}
where we have used identities (\ref{iden-1}), (\ref{iden-3}) and (\ref{iden-2}). Integrating (\ref{L2-1}) over the interval $[0,t]$ yields
\begin{align}  \label{L2-4}
\|w(t)\|_2^2 + \|\mathbf u(t)\|_2^2 +\int_0^t \left(\frac{2}{Re} \left(\|\nabla_h w\|_2^2 + \|\nabla_h \mathbf u\|_2^2 \right)  + 2\epsilon^2 \|\psi_z\|_2^2 \right)ds  =   \|w_0\|_2^2 + \|\mathbf u_0\|_2^2.
\end{align}

\vspace{0.1 in}

\subsection{Estimate for $\|\omega\|_2^2$} Taking the inner product of (\ref{model-2}) with $\omega$ yields
\begin{align}  \label{omega-1}
\frac{1}{2}\frac{d}{dt} \|\omega\|_2^2 
+ \frac{1}{Re}\|\nabla_h \omega\|_2^2 + \epsilon^2 \|\mathbf u_z\|_2^2 = (w_z,\omega),
\end{align}
where (\ref{iden-1}) and (\ref{iden-3}) have been used. Thanks to (\ref{iden-4}), we have
\begin{align}   \label{omega-2}
&(w_z,\omega)=\int_{\Omega} w_z (-\Delta_h \psi) dx dy dz = -\int_{\Omega} \nabla_{h}w \cdot  \nabla_h \psi_z  dx dy dz  \notag\\
&\leq \|\nabla_h w\|_2\|\nabla_h \psi_z\|_2 = \|\nabla_h w\|_2 \|\mathbf u_z\|_2 \leq \frac{\epsilon^2}{2} \|\mathbf u_z\|_2^2 + \frac{1}{2\epsilon^2}  \|\nabla_h w\|_2^2. 
\end{align}
Combining (\ref{omega-1}) and (\ref{omega-2}) implies
\begin{align}   \label{omega-3}
\frac{d}{dt} \|\omega\|_2^2 
+ \frac{2}{Re}\|\nabla_h \omega\|_2^2 + \epsilon^2 \|\mathbf u_z\|_2^2 
\leq \frac{1}{\epsilon^2} \|\nabla_h w\|_2^2.
\end{align}
By integrating (\ref{omega-3}) over the interval $[0,t]$, we obtain
\begin{align}  \label{omega}
&\|\omega(t)\|_2^2 + \int_0^t \left(\frac{2}{Re}\|\nabla_h \omega\|_2^2 + \epsilon^2 \|\mathbf u_z\|_2^2  \right) ds
 \leq \|\omega_0\|_2^2 + \frac{1}{\epsilon^2} \int_0^t  \|\nabla_h w\|_2^2  ds \notag\\
&\leq \|\omega_0\|_2^2 + \frac{Re}{2\epsilon^2} \left(\|\omega_0\|_2^2 + \|\mathbf u_0\|_2^2 \right),
 \end{align}
where the last inequality is due to (\ref{L2-4}).

\vspace{0.1 in}

\subsection{An anisotropic Ladyzhenskaya type inequality} 
We state here the following anisotropic Ladyzhenskaya type inequality which will be useful in subsequent {\it a priori} estimates. It is worth mentioning that similar inequalities can be found in \cite{CT11}. However, for the sake of completeness we present the proof of this technical lemma in the appendix. 
\begin{lemma}   \label{lemma1}
Let $f\in H^1(\Omega)$, $g\in H^1_h(\Omega)$ and $h\in L^2(\Omega)$. Then
\begin{align*}
\int_{\Omega} |fgh| dx dy dz  \leq C \left(\|f\|_2 + \|\nabla_h f\|_2 \right)^{\frac{1}{2}}  \left(\|f\|_2 + \|f_z\|_2 \right)^{\frac{1}{2}} 
 \|g\|_2^{\frac{1}{2}} \left(\|g\|_2  + \|\nabla_h g\|_2 \right)^{\frac{1}{2}} \|h\|_2.
\end{align*}
\end{lemma}

\vspace{0.1 in}

\subsection{Estimate for $\|\nabla_h w\|_2$} 
Taking the inner product of (\ref{model-1}) with $-\Delta_h w$ yields
\begin{align*} 
&\frac{1}{2} \frac{d}{dt}\|\nabla_h w\|_2^2 +\frac{1}{Re} \|\Delta_h w\|_2^2    \notag\\
&\leq \int_{\Omega} |\mathbf u \cdot \nabla_h w \Delta_h w|  dx dy dz 
+ \|\psi_z\|_2 \|\Delta_h w\|_2  \notag\\
&\leq C   \|\omega\|_2^{1/2} (\|\mathbf u\|_2 + \|\mathbf u_z\|_2)^{1/2}   \|\nabla_h w\|_2^{1/2}  \|\Delta_h w\|_2^{3/2} +\|\psi_z\|_2 \|\Delta_h w\|_2,
\end{align*}
where we have used Lemma \ref{lemma1} and the Poincar\'e inequality since $\int_{[0,L]^2}  \textbf u dx dy =  \int_{[0,L]^2} (\psi_y, - \psi_x)^{tr} dx dy=0$ 
and $\int_{[0,L]^2} \nabla_h w dx dy =0$. 

By employing the Young's inequality, we obtain
\begin{align*}   
\frac{d}{dt}\|\nabla_h w\|_2^2 +\frac{1}{Re} \|\Delta_h w\|_2^2  \leq C\|\omega\|_2^2(\|\mathbf u\|_2^2 + \|\mathbf u_z\|_2^2) \|\nabla_h w\|_2^2 + C \|\psi_z\|_2^2. 
\end{align*}
Thanks to the Gronwall's inequality, we have
\begin{align}  \label{H1w}
\|\nabla_h w(t)\|_2^2 + \frac{1}{Re}\int_0^t \|\Delta_h w\|_2^2 ds 
&\leq  C\left(\|\nabla_h w_0\|_2^2 + \int_0^t \|\psi_z\|_2^2  ds\right)
e^{\int_0^t C \|\omega\|_2^2 (\|\mathbf u\|_2^2 + \|\mathbf u_z\|_2^2) ds} \notag\\
&\leq C\left(\|w_0\|_2, \|\nabla_h w_0\|_2, \|\omega_0\|_2 \right).
\end{align}
The uniform bound (\ref{H1w}) is due to estimates (\ref{L2-4}) and (\ref{omega}).

\vspace{0.1 in}

\subsection{Estimate for $\|w_z\|_2^2 + \|\mathbf u_z\|_2^2$} 
We take the $L^2(\Omega)$ inner product of (\ref{model-1})-(\ref{model-2}) with $(-w_{zz}, -\psi_{zz})^{tr}$. After conducting integration by parts, one has
\begin{align}   \label{H1z-1}
&\frac{1}{2}\frac{d}{dt}\left(\|w_z\|_2^2 + \|\mathbf u_z\|_2^2\right) + \frac{1}{Re} \left(\|\nabla_h w_z\|_2^2 + \|\omega_z\|_2^2 \right) + \epsilon^2 \|\psi_{zz}\|_2^2   \notag\\
&\leq \int_{\Omega} |{\mathbf u}_z \cdot \nabla_h w w_z| dx dy dz +  
\int_{\Omega} |{\mathbf u}_z \cdot \nabla_h \psi_z \omega| dx dy dz
+   \int_{\Omega} |\mathbf u \cdot \nabla_h \psi_z \omega_z| dx dy dz.
\end{align}
Next, we estimate each term on the right-hand side of (\ref{H1z-1}).

By Lemma \ref{lemma1} with $f=\nabla_h w$, $g=\mathbf u_z$ and $h=w_z$, and along with the Poincar\'e inequality, we obtain
\begin{align}   \label{H1z-2}
&\int_{\Omega} |{\mathbf u}_z \cdot \nabla_h w w_z| dx dy dz  \notag\\
&\leq C \|\Delta_h w\|_2^{1/2}    \left(\|\nabla_h w\|_2 + \|\nabla_h w_z\|_2 \right)^{1/2}   
 \|\mathbf u_z\|_2^{1/2}   \|\omega_z\|_2  ^{1/2} \|w_z\|_2 \notag\\
&\leq \frac{1}{6Re}\left(\|\nabla_h w_z\|_2^2 + \|\omega_z\|_2^2 \right)
+ C \left(\|\Delta_h w\|_2^2 + \|\mathbf u_z\|_2^2\right) \left(\|w_z\|_2^2+1\right).
\end{align}
Also using Lemma \ref{lemma1} with $f=\omega$, $g=\mathbf u_z$ and $h=\nabla_h \psi_z$ gives us
\begin{align}  \label{H1z-3}
&\int_{\Omega} |{\mathbf u}_z \cdot \nabla_h \psi_z \omega| dx dy dz   \notag\\
&\leq C \|\nabla_h \omega\|_2^{1/2} \left(\|\omega\|_2 +  \|\omega_z\|_2\right)^{1/2} \|\mathbf u_z\|_2^{3/2} \|\omega_z\|_2^{1/2} \notag\\
&\leq \frac{1}{6Re} \|\omega_z\|_2^2 + C \left(\|\nabla_h \omega\|_2^2 + \|\mathbf u_z\|_2^2+1  \right) \|\mathbf u_z\|_2^2.
\end{align}
In addition, due to Lemma \ref{lemma1} with $f=\mathbf u$, $g=\nabla_h \psi_z$ and $h=\omega_z$, one has
\begin{align} \label{H1z-4}
&\int_{\Omega} |\mathbf u \cdot \nabla_h \psi_z \omega_z| dx dy dz  \notag\\
&\leq C \|\omega\|_2^{1/2} \left(\|\mathbf u\|_2 + \|\mathbf u_z\|_2\right)^{1/2} 
\|\mathbf u_z\|_2^{1/2} \|\omega_z\|_2^{3/2} \notag\\
&\leq \frac{1}{6Re}\|\omega_z\|_2^2 + C \|\omega\|_2^2  \left(\|\mathbf u\|_2^2 + \|\mathbf u_z\|_2^2 \right)   \|\mathbf u_z\|_2^2.
\end{align}

Apply estimates (\ref{H1z-2})-(\ref{H1z-4}) to the inequality (\ref{H1z-1}) yields
\begin{align*}
&\frac{d}{dt}\left(\|w_z\|_2^2 + \|\mathbf u_z\|_2^2 \right) + 
\frac{1}{Re} \left(\|\nabla_h w_z\|_2^2 + \|\omega_z\|_2^2 \right) 
+ \epsilon^2 \|\psi_{zz}\|_2^2    \notag\\
&\leq C \left(\|\Delta_h w\|_2^2 + \|\mathbf u_z\|_2^2\right) \left(\|w_z\|_2^2+1\right)    \notag\\
&\hspace{0.2 in}+C \left(\|\nabla_h \omega\|_2^2 + \|\mathbf u_z\|_2^2   + \|\omega\|_2^2\|\mathbf u\|_2^2 +   \|\omega\|_2^2 \|\mathbf u_z\|_2^2 \right)   \|\mathbf u_z\|_2^2.
\end{align*}
Thanks to Gronwall's inequality, we obtain
\begin{align}   \label{H1z-11}
&\|w_z(t)\|_2^2 + \|\mathbf u_z(t)\|_2^2 +\int_0^t  \left[\frac{1}{Re} \left(\|\nabla_h w_z\|_2^2 + \|\omega_z\|_2^2 \right)  + \epsilon^2 \|\psi_{zz}\|_2^2 \right]ds \notag\\
&\leq \left( \|\partial_z{w_0}\|_2^2 + \|\partial_z \mathbf u_0\|_2^2  + C\int_0^t \left(\|\Delta_h w\|_2^2 + \|\mathbf u_z\|_2^2  
+\|\omega\|_2^2\|\mathbf u\|_2^2 \|\mathbf u_z\|_2^2   \right) ds\right)   \notag\\
&\hspace{0.2 in} \exp\left\{C\int_0^t  \big(\|\Delta_h w\|_2^2 + \|\mathbf u_z\|_2^2 + \|\nabla_h \omega\|_2^2
+ \|\omega\|_2^2 \|\mathbf u_z\|_2^2  \big) ds\right\} \notag\\
&\leq C(\|w_0\|_{H^1}, \|\mathbf u_0\|_{H^1}).
\end{align}
The uniform bound (\ref{H1z-11}) is due to (\ref{L2-4}), (\ref{omega}) and (\ref{H1w}).

\vspace{0.2 in}

\section{Rigorous justification of the {\it a priori} estimates and the existence of strong solutions}      \label{sec-existence}
This section is devoted to proving the existence of global strong solutions for the model (\ref{model-1})-(\ref{model-4}) by assuming the initial data 
$(\mathbf u_0,w_0)^{tr} \in (H^1(\Omega))^3$. 
We employ the standard Galerkin method and use the analogue of the {\it a priori} estimates that were established in section \ref{sec-est}.

Let $e_{\mathbf j} = \exp\left(2\pi i [(j_1 x+j_2 y)/L + j_3 z] \right)$ for $\mathbf j=(j_1,j_2,j_3)^{tr}$. For $m\in \mathbb N$, let $P_m\left(L^2(\Omega)\right)$ be a subspace of $L^2(\Omega)$ spanned by $\{e_{\mathbf j}\}_{|\mathbf j|\leq m}$. Also, for any $L^2(\Omega)$ function 
$f=\sum_{\mathbf j\in \mathbb Z^3} \alpha_{\mathbf j} e_{\mathbf j}$, with $\alpha_{\mathbf j}=(f,e_{\mathbf j})$, we write 
$P_m f=\sum_{|\mathbf j|\leq m} \alpha_{\mathbf j} e_{\mathbf j}.$

Let us consider the Galerkin approximation for our model (\ref{model-1})-(\ref{model-4}):
\begin{align}     
&\partial_t w_m + P_m\left(\mathbf u_m \cdot \nabla_h w_m\right) - \partial_z \psi_m = \frac{1}{Re} \Delta_h w_m, \label{app-1} \\
&\partial_t \omega_m + P_m\left(\mathbf u_m \cdot \nabla_h \omega_m\right) - \partial_z w_m = \frac{1}{Re} \Delta_h \omega_m  + \epsilon^2 \partial_{zz} \psi_m , \label{app-2}\\
&\nabla_h \cdot \mathbf u_m = 0,   \label{app-3}   \\
& \mathbf u_m(0)= P_m \mathbf u_0,    \;\;\; w_m(0)=P_m w_0,    \label{app-4}
\end{align}  
where $\mathbf u_m, w_m \in P_m\left(L^2(\Omega)\right)$ and 
$\omega_m = \nabla_h \times \mathbf u_m$, $\psi_m = (-\Delta_h)^{-1} \omega_m$ with $\int_{[0,L]^2} \psi_m dx dy =0$.

For each $m\geq 1$, the Galerkin approximation (\ref{app-1})-(\ref{app-4}) corresponds to a first order system of ordinary differential equations with quadratic nonlinearity. Therefore, by the theory of ordinary differential equations, there exists some $T_m>0$ such that system (\ref{app-1})-(\ref{app-4}) admits a unique solution $(\mathbf u_m, w_m)^{tr}$ on $[0,T_m]$. Since $\mathbf u_m$ and $w_m$ have finitely many modes, they are smooth functions, and therefore all of the {\it a priori} estimates established in section \ref{sec-est} are valid for the Galerkin approximate solution $(\mathbf u_m, w_m)^{tr}$. In particular, the $H^1(\Omega)$ norm of $(\mathbf u_m, w_m)^{tr}$ is uniformly bounded for all time.  
Hence, the Galerkin approximate solution $(\mathbf u_m, w_m)^{tr}$ exists globally in time, in particular, over $[0,T]$, for every $T>0$.

Furthermore, by the {\it a priori} estimates in section \ref{sec-est}, one has the following uniform bounds for the sequence of the Galerkin approximate solutions.
\begin{align}    
&\mathbf u_m, \; w_m  \;\;\text{are uniformly bounded in}\;\;  L^{\infty}(0,T; H^1(\Omega));   \label{Gal-2} \\
&\nabla_h \omega_m,\; \Delta_h w_m,\; \partial_z{\omega_m}, \; \nabla_h \partial_z {w_m}, \; \partial_{zz}\psi_m \;\text{are uniformly bounded in}\; L^2(\Omega \times (0,T))  \label{Gal-3}.
\end{align} 
Therefore, there exist a subsequence, denoted also by $\mathbf u_m$, $w_m$, $\omega_m$, $\psi_m$, and corresponding limits, $\mathbf u$, $w$, $\omega$, and $\psi$, respectively, such that 
\begin{align}    
&\mathbf u_m \rightarrow \mathbf u, \;\; w_m\rightarrow w,  \;\; \text{weakly$^*$ in} \;\; L^{\infty}(0,T;H^1(\Omega));   \label{Gal-5'}  \\
&\nabla_h \omega_m \rightarrow \nabla_h \omega,       \;\; 
\Delta_h w_m \rightarrow \Delta_h w,     \;\; \text{weakly in}  \;\; L^2(\Omega \times (0,T));    \label{Gal-5''}  \\
&\partial_z{\omega_m} \rightarrow \partial_z \omega, \;\; \nabla_h \partial_z {w_m} \rightarrow \nabla_h w_z,  
\;\; \partial_{zz}\psi_m \rightarrow \partial_{zz} \psi, \;\; \text{weakly in} \;\;   L^2(\Omega \times (0,T)).   \label{Gal-5'''}
\end{align}

Moreover, due to the {\it a priori} estimates in section \ref{sec-est}, we find that 
\begin{align}   \label{unibnd}
\sup_{0\leq t \leq T} \left(\|\mathbf u_m(t)\|_{H^1(\Omega)}^2 + \|w_m(t)\|_{H^1(\Omega)}^2\right) \leq K
\end{align}
where $K$ is independent of $T$, but depends only on parameters $Re$, $\epsilon$, $L$ as well as the $H^1$-norm, $\|\mathbf u_0\|_{H^1(\Omega)}$ and $\|w_0\|_{H^1(\Omega)}$ of the initial data.
Also thanks to the $\text{weak-}*$ convergence stated in (\ref{Gal-5'}), one has $\|\mathbf u\|_{L^{\infty}(0,T;H^1(\Omega))} \leq \liminf_{m\rightarrow \infty} \|\mathbf u_m\|_{L^{\infty}(0,T;H^1(\Omega))}$ and 
$\|w\|_{L^{\infty}(0,T;H^1(\Omega))} \leq \liminf_{m\rightarrow \infty} \|w_m\|_{L^{\infty}(0,T;H^1(\Omega))}$. Therefore, we obtain from (\ref{unibnd}) that
$$\sup_{0\leq t \leq T} \left(\|\mathbf u(t)\|_{H^1(\Omega)}^2 + \|w(t)\|_{H^1(\Omega)}^2\right) \leq K.$$

In order to obtain the strong convergence of the approximate solutions, we shall derive uniform bounds for $\partial_t w_m$ and $\partial_t \mathbf u_m$. First, we claim that the sequence $\partial_t w_m$ is uniformly bounded in $L^2(\Omega \times (0,T))$. Indeed, for any function $\varphi\in L^{4/3}(0,T;L^2(\Omega))$, we use Lemma \ref{lemma1} to estimate
\begin{align}  \label{Gal-5}
&\int_0^T \int_{\Omega} |(\mathbf u_m \cdot \nabla_h w_m) \varphi | dx dy dz  dt    \notag\\
&\leq C\int_0^T \|\omega_m\|_2^{1/2} \left(\|\mathbf u_m\|_2+ \|\partial_z \mathbf u_m\|_2\right)^{1/2} \|\nabla_h w_m\|_2^{1/2}  
\|\Delta_h w_m\|_2^{1/2} \|\varphi\|_2 dt  \notag\\
&\leq C\sup_{t\in [0,T]}  \left(\|\omega_m\|_2^{1/2} \left(\|\mathbf u_m\|_2+ \|\partial_z \mathbf u_m\|_2\right)^{1/2} \|\nabla_h w_m\|_2^{1/2} \right) \notag\\
&\hspace{0.5 in} \cdot \left(\int_0^T  \|\Delta_h w_m\|_2^2  dt\right)^{1/4} \left(\int_0^T \|\varphi\|_2^{4/3}dt\right)^{3/4} \notag\\
& \leq C(\|\mathbf u_0\|_{H^1},  \|w_0\|_{H^1})  \|\varphi\|_{L^{4/3}(0,T;L^2(\Omega))},
\end{align}
where the last inequality is due to the {\it a priori} estimates (\ref{L2-4}), (\ref{omega}), (\ref{H1w}) and (\ref{H1z-11}). 
Consequently, the sequence
\begin{align}    \label{Gal-6'}
\mathbf u_m \cdot \nabla_h w_m \;\;\text{is uniformly bounded in}\;\; L^4(0,T;L^2(\Omega)). 
\end{align}
As a result, from (\ref{Gal-2})-(\ref{Gal-3}) and (\ref{Gal-6'}), we obtain from (\ref{app-1}) that the sequence
\begin{align}    \label{Gal-6}
\partial_t w_m \;\;\text{is uniformly bounded in}\;\; L^2(\Omega \times (0,T)). 
\end{align}

Next, we show that $\partial_t \mathbf u_m $ is uniformly bounded in $L^2(\Omega \times (0,T))$. Recall the Hilbert space 
$H^1_h(\Omega)=\{f\in L^2(\Omega):     \nabla_h f \in L^2(\Omega)\}$ associated with the norm $\|f\|^2_{H^1_h(\Omega)}=\|f\|_2^2+\|\nabla_h f\|_2^2$. 
For any function $\phi \in L^2(0,T;H^1_h(\Omega))$, we apply Lemma \ref{lemma1} in order to estimate
\begin{align}     \label{Gal-7}
&\int_0^T \int_{\Omega} |(\mathbf u_m \cdot \nabla_h \omega_m) \phi| dx dy dz dt \notag\\
&\leq  C\int_0^T \|\omega_m\|_2^{1/2} \left(\|\mathbf u_m\|_2 + \|\partial_z \mathbf u_m\|_2\right)^{1/2} \|\nabla_h \omega_m\|_2 
\|  \phi   \|_2^{1/2}    \left(\|   \phi    \|_2 + \|\nabla_h    \phi     \|_2\right)^{1/2}   dt \notag\\
&\leq C\sup_{t\in[0,T]}  \left(\|\omega_m\|_2^{1/2} \left(\|\mathbf u_m\|_2 + \|\partial_z\mathbf u_m\|_2\right)^{1/2} \right)  \notag\\
&\hspace{0.5 in}  \cdot \left(\int_0^T \|\nabla_h \omega_m\|_2^2 dt \right)^{1/2} \left(\int_0^T \left(\|   \phi    \|_2^2 + \|\nabla_h   \phi     \|_2^2 \right) dt\right)^{1/2} \notag\\
&\leq C(\|\mathbf u_0\|_{H^1},  \|w_0\|_{H^1})  \|\phi\|_{L^2(0,T;H^1_h(\Omega))},
\end{align}
where we have used the {\it a priori} estimates (\ref{L2-4}), (\ref{omega}) and (\ref{H1z-11}). Therefore, the sequence 
\begin{align}  \label{Gal-7'}
\mathbf u_m \cdot \nabla_h \omega_m \;\;\text{is uniformly bounded in}\;\; L^2(0,T;H^1_h(\Omega)'),
\end{align}
where $(H^1_h(\Omega))'$ is the dual space of $H^1_h(\Omega)$.
Consequently, according to (\ref{Gal-3}) and (\ref{Gal-7'}), we obtain from (\ref{app-2}) that the sequence
\begin{align}  \label{Gal-8'}
\partial_t \omega_m  \;\;\text{is uniformly bounded in}\;\;  L^2(0,T;H^1_h(\Omega)'),
\end{align}
and thus
\begin{align}  \label{Gal-8}
\partial_t \mathbf u_m  \;\;\text{is uniformly bounded in}\;\;  L^2(\Omega \times (0,T)).
\end{align}

Then, we infer from (\ref{Gal-6}) and (\ref{Gal-8}) that there is a subsequence such that
\begin{align}  \label{Gal-111}
\partial_t w_m \rightarrow \partial_t w, \;\;  \partial_t \mathbf u_m \rightarrow \partial_t \mathbf u \;\;\; \text{weakly in} \;\; L^2(\Omega \times (0,T)).  
\end{align}

By (\ref{Gal-2}), (\ref{Gal-6}), (\ref{Gal-8}), and thanks to the Aubin's compactness theorem, we have, for a subsequence, the following strong convergence holds:
\begin{align}     \label{Gal-12}
&\mathbf u_m \rightarrow \mathbf u, \;\; w_m\rightarrow w  \;\; \text{in} \;\; L^2(\Omega \times (0,T)).
\end{align}

Next, we show the convergence of the nonlinear terms in (\ref{app-1})-(\ref{app-2}). Let $\eta$ be a trigonometric polynomial with continuous coefficients. For $m$ larger than the degree of $\eta$ we have
\begin{align}    \label{Gal-13}
&\int_0^T \int_{\Omega} P_m(\mathbf u_m \cdot \nabla_h \omega_m)\eta dx dy dz   \notag\\
&= \int_0^T \int_{\Omega} (\mathbf u \cdot \nabla_h \omega_m) \eta dx dy dz + \int_0^T \int_{\Omega} \left((\mathbf u_m-\mathbf u) \cdot \nabla_h \omega_m\right) \eta dx dy dz.
\end{align}
Since $\nabla_h \omega_m\rightarrow \nabla_h \omega$ weakly in $L^2(\Omega \times (0,T))$, $\mathbf u_m \rightarrow \mathbf u$ in $L^2(\Omega \times (0,T))$,
and $\nabla_h \omega_m$ is uniformly bounded in $L^2(\Omega \times (0,T))$, we can pass to the limit in (\ref{Gal-13}):
\begin{align}     \label{Gal-14}
\lim_{m\rightarrow \infty}  \int_0^T \int_{\Omega} P_m (\mathbf u_m \cdot \nabla_h \omega_m) \eta dx dy dz
=  \int_0^T \int_{\Omega} (\mathbf u \cdot \nabla_h \omega)\eta dx dy dz.
\end{align}
An analogous argument yields
\begin{align}    \label{Gal-15}
&\lim_{m\rightarrow \infty}  \int_0^T \int_{\Omega} P_m (\mathbf u_m \cdot \nabla_h w_m) \eta dx dy dz
=  \int_0^T \int_{\Omega} (\mathbf u \cdot \nabla_h w) \eta dx dy dz. 
\end{align}

Therefore, due to (\ref{Gal-5'})-(\ref{Gal-5'''}), (\ref{Gal-111}), (\ref{Gal-14}) and (\ref{Gal-15}), we pass to the limit for the Galerkin approximate equations (\ref{app-1})-(\ref{app-3}). It follows that
\begin{align}    
&\int_0^T \int_{\Omega} \left(\partial_t w + \mathbf u \cdot \nabla_h w - \partial_z \psi  - \frac{1}{Re} \Delta_h w \right) \eta dx dy dz dt=0,  \label{Gal-22}   \\
&\int_0^T \int_{\Omega} \left(\partial_t \omega + \mathbf u \cdot \nabla_h \omega - \partial_z w -\frac{1}{Re} \Delta_h \omega  - \epsilon^2 \partial_{zz} \psi \right)\eta dx dy dz dt=0,    \label{Gal-221} \\
&\int_0^T \int_{\Omega}  (\nabla_h \cdot u) \eta dx dy dz dt =0,   \label{Gal-222}
\end{align}    
for any trigonometric polynomial $\eta$ with continuous coefficients.

By applying Lemma \ref{lemma1} as the arguments in (\ref{Gal-5}), we can deduce that $\mathbf u \cdot \nabla_h w\in L^4(0,T;L^2(\Omega))$. Then, since $\partial_t w$, $\partial_{zz}\psi$ and $\Delta_h w \in L^2(\Omega \times (0,T))$, one has
\begin{align}  \label{Gal-22'}
\partial_t w + \mathbf u \cdot \nabla_h w - \partial_z \psi - \frac{1}{Re} \Delta_h w\in L^2(\Omega \times (0,T)).
\end{align}
Also, using Lemma \ref{lemma1} as the estimates in (\ref{Gal-7}), one may derive that $\mathbf u \cdot \nabla_h \omega \in L^2(0,T; H^1_h(\Omega)')$. 
Since $\partial_z w$, $\Delta_h \omega$, $\partial_{zz}\psi\in L^2(\Omega \times (0,T))$ and $\partial_t \omega\in L^2(0,T; H^1_h(\Omega)') $, we obtain
\begin{align}      \label{Gal-22''}
\partial_t \omega + \mathbf u \cdot \nabla_h \omega - \partial_z w -\frac{1}{Re} \Delta_h \omega  - \epsilon^2 \partial_{zz} \psi \in L^2(0,T;H^1_h(\Omega)').
\end{align}

On account of (\ref{Gal-22'}) and (\ref{Gal-22''}), we obtain from (\ref{Gal-22})-(\ref{Gal-222}) that
\begin{align}     \label{Gal-23}
&\partial_t w + \mathbf u \cdot \nabla_h w - \partial_z \psi =  \frac{1}{Re} \Delta_h w,     \;\;\text{in}\;\; L^2(\Omega \times (0,T));\\
&\partial_t \omega + \mathbf u \cdot \nabla_h \omega - \partial_z w = \frac{1}{Re} \Delta_h \omega  + \epsilon^2 \partial_{zz} \psi ,
\;\;\text{in}\;\; L^2(0,T; H^1_h(\Omega)'),  \label{Gal-231}
\end{align}    
with $\nabla_h \cdot \mathbf u=0$.

It follows from (\ref{Gal-23})-(\ref{Gal-231}) that 
\begin{align}   \label{Gal-1111}
\mathbf u, w \in C(0,T;L^2(\Omega)).
\end{align}
Due to (\ref{Gal-1111}) and (\ref{Gal-12}), one has, for every $t\in [0,T]$,
$\mathbf u_m(t) \rightarrow \mathbf u(t)$ and $w_m(t)\rightarrow w(t)$ in $L^2(\Omega)$. 
In particular, $\mathbf u_m(0) \rightarrow \mathbf u(0)$ and $w_m(0)\rightarrow w(0)$. 
On the other hand, by (\ref{app-4}), we find that $\mathbf u_m(0) \rightarrow \mathbf u_0$ and $w_m(0)\rightarrow w_0$. As a result, $(\mathbf u,w)^{tr}$ satisfies the desired initial condition: $\mathbf u(0)=\mathbf u_0$ and $w(0)=w_0$.

Finally, due to the regularity of solutions, we can multiply (\ref{Gal-23})-(\ref{Gal-231}) by $(w,\psi)^{tr}$ and integrate the result over $\Omega \times (0,t)$ for $t\in [0,T]$. Then the energy identity (\ref{energy}) follows.

\vspace{0.2 in}

\section{Uniqueness of strong solutions}  \label{sec-unique}
This section is devoted to proving that strong solutions for the system (\ref{model-1})-(\ref{model-4}) are unique and depend continuously on the initial data. Assume there are two strong solutions $(\mathbf u_1,w_1)^{tr}$ and $(\mathbf u_2,w_2)^{tr}$ on $[0,T]$ in the sense of Definition \ref{def-sol}. 
Set $\mathbf u=\mathbf u_1-\mathbf u_2$ and $w=w_1-w_2$. 
Therefore, 
\begin{align}     
&\frac{\partial w}{\partial t} + \mathbf u \cdot \nabla_h w_1+ \mathbf u_2 \cdot \nabla_h w -  \frac{\partial \psi}{\partial z}  = \frac{1}{Re} \Delta_h w,  \;\;  \text{in}\; L^2(\Omega \times (0,T)); \label{unique}    \\
&\frac{\partial \omega}{\partial t} + \mathbf u \cdot \nabla_h \omega_1 + \mathbf u_2 \cdot \nabla_h \omega 
-\frac{\partial w}{\partial z} = \frac{1}{Re} \Delta_h \omega  + \epsilon^2    \frac{\partial^2 \psi}{\partial z^2}  , 
\;\;\text{in}\; L^2(0,T;H^1_h(\Omega)'), \label{unique0}
\end{align}
with $\nabla_h \cdot \mathbf u=0$.

Since $\mathbf u$ and $w$ satisfy the regularity (\ref{main-2}), we can multiply (\ref{unique})-(\ref{unique0}) by $(w,\psi)^{tr}$ and integrate over $\Omega$. By using (\ref{iden-0}), (\ref{iden-1}), (\ref{iden-3}) and (\ref{iden-2}), we obtain, for a.e. $t\in [0,T]$,
\begin{align}   \label{unique-1}
&\frac{1}{2}\frac{d}{dt}\left(\|w\|_2^2 + \|\mathbf u\|_2^2\right) 
+ \frac{1}{Re}\left(\|\nabla_h w\|_2^2 + \|\nabla_h \mathbf u\|_2^2\right) + \epsilon^2\|\psi_z\|_2^2 \notag\\
&\leq \int_{\Omega} |(\mathbf u \cdot \nabla_h w) w_1| dx dy dz+\int_{\Omega} |(\mathbf u_2 \cdot \nabla_h \psi) \omega| dx dy dz.
\end{align}
Next we estimate the two integrals on the right-hand side of (\ref{unique-1}).

Using Lemma 1 with $f=w_1$, $g=\mathbf u$ and $h=\nabla_h w$, we obtain
\begin{align}    \label{unique-2}
&\int_{\Omega} |(\mathbf u \cdot \nabla_h w) w_1| dx dy dz   \notag\\
&\leq C \left(\|w_1\|_2 + \|\nabla_h w_1\|_2 \right)^{1/2}  \left(\|w_1\|_2 + \|\partial_z w_1\|_2\right)^{1/2}
\|\mathbf u\|_2^{1/2} \|\nabla_h \mathbf u\|_2^{1/2} \|\nabla_h w\|_2 \notag\\
&\leq \frac{1}{4Re} \left(\|\nabla_h w\|_2^2 + \|\nabla_h \mathbf u\|_2^2\right)+
C\left(\|w_1\|_2^2 + \|\nabla_h w_1\|_2^2 \right)  \left(\|w_1\|_2^2 + \|\partial_z w_1\|_2^2\right) \|\mathbf u\|_2^2.
\end{align}
Also, using Lemma 1 with $f=\mathbf u_2$, $g=\nabla_h \psi$, $h=\omega$, we have
\begin{align}  \label{unique-3}
\int_{\Omega} |(\mathbf u_2 \cdot \nabla_h \psi) \omega| dx dy dz 
&\leq C    \|\omega_2\|_2^{1/2}        (\|\mathbf u_2\|_2+ \|\partial_z \mathbf u_2\|_2)^{1/2} \|\mathbf u\|_2^{1/2} \|\nabla_h \mathbf u\|_2^{3/2}  \notag\\
&\leq  \frac{1}{4Re} \|\nabla_h \mathbf u\|_2^2 +   C    \|\omega_2\|_2^2      (\|\mathbf u_2\|_2^2+ \|\partial_z \mathbf u_2\|_2^2)\|\mathbf u\|_2^2.
\end{align}

Now, we combine the estimates (\ref{unique-1})-(\ref{unique-3}) to deduce, for a.e. $t\in [0,T]$,
\begin{align*} 
&\frac{d}{dt}\left(\|w\|_2^2 + \|\mathbf u\|_2^2  \right) 
+ \frac{1}{Re}\left(\|\nabla_h w\|_2^2 + \|\nabla_h \mathbf u\|_2^2\right) + \epsilon^2\|\psi_z\|_2^2    \notag\\
&\leq  C  \Big[\left(\|w_1\|_2^2 + \|\nabla_h w_1\|_2^2 \right)  \left(\|w_1\|_2^2 + \|\partial_z w_1\|_2^2\right)+\|\omega_2\|_2^2 (\|\mathbf u_2\|_2^2+ \|\partial_z \mathbf u_2\|_2^2)\Big]\|\mathbf u\|_2^2.
\end{align*}
By Gronwall's inequality, it follows that
\begin{align}   \label{unique-4}
&\|w(t)\|_2^2 + \|\mathbf u(t)\|_2^2    \notag\\
&\leq   \left(\|w(0)\|_2^2 + \|\mathbf u(0)\|_2^2\right) 
e^{C\int_0^t \left(\|w_1\|_2^2 + \|\nabla_h w_1\|_2^2 \right)  \left(\|w_1\|_2^2 + \|\partial_z w_1\|_2^2\right)+\|\omega_2\|_2^2 (\|\mathbf u_2\|_2^2+ \|\partial_z \mathbf u_2\|_2^2)ds} \notag\\
&\leq    \left(\|w(0)\|_2^2 + \|\mathbf u(0)\|_2^2\right) e^{t C\left(\|w_1(0)\|_{H^1},\|\mathbf u_2(0)\|_{H^1}\right)},
\end{align}
for any $t\in [0,T]$. In particular, if $(\mathbf u(0),w(0))^{tr}=0$, i.e., the initial values of the two solutions $(\mathbf u_1,w_1)^{tr}$ and $(\mathbf u_2,w_2)^{tr}$ coincide, then (\ref{unique-4}) implies $\|w(t)\|_2^2 + \|\mathbf u(t)\|_2^2=0$ for all $t\in [0,T]$. This completes the proof for the uniqueness of strong solutions.

To see the continuous dependence on the initial data, we let  $({\tilde{\mathbf u}}_0, {{\tilde w}}_0)^{tr} \in (H^1(\Omega))^3$ and take a bounded sequence $\{(\mathbf u_0^n,w_0^n)^{tr}\}$ of initial data in $H^1(\Omega)$ such that $(\mathbf u_0^n,w_0^n)^{tr} \rightarrow ({\tilde{\mathbf u}}_0,\tilde w_0)^{tr}$ in $L^2(\Omega)$, and 
$\|\mathbf u_0^n\|_{H^1}$, $\|w_0^n\|_{H^1}$, $\|{\tilde{\mathbf u}}_0\|_{H^1}$, $\|\tilde w_0\|_{H^1} \leq M$ for some $M>0$. Denote the corresponding strong solutions by $(\mathbf u^n,w^n)^{tr}$ and $({\tilde{\mathbf u}},\tilde w)^{tr}$, respectively. 
Then, on account of (\ref{unique-4}), we have, for all $t\in [0,T]$,
\begin{align*}
\|\tilde w - w^n\|_2^2 + \|{\tilde{\mathbf u}}-\mathbf u^n\|_2^2 
&\leq  \left(\|\tilde w_0 - w_0^n  \|_2^2 + \|{\tilde{\mathbf u}}_0  -\mathbf u^n_0  \|_2^2\right) e^{t C\left(\|\tilde w_0\|_{H^1},\|\mathbf u^n_0\|_{H^1}\right)} \\
&\leq  \left(\|\tilde w_0 - w_0^n  \|_2^2 + \|{\tilde{\mathbf u}}_0  -\mathbf u^n_0  \|_2^2\right) e^{T\cdot C(M)}.
\end{align*}
It follows that $(\mathbf u^n,w^n)^{tr} \rightarrow ({\tilde{\mathbf u}},\tilde w)^{tr}$ in $C([0,T];L^2(\Omega))$. This completes the proof for the continuous dependence on the initial data with respect to the $L^2$-norm for the strong solutions.

\vspace{0.2 in}

\section{Appendix}
We prove the anisotropic Ladyzhenskaya type inequality stated in Lemma \ref{lemma1}.
\begin{proof}
It suffices to prove the inequality in Lemma \ref{lemma1} for smooth periodic functions, and then pass to the limit using a standard density argument. 
Recall $\Omega=[0,L]^2 \times [0,1]$. For a fixed $(x,y)\in [0,L]^2$ and for every $z$, $\sigma \in [0,1]$, we have
\begin{align}  \label{lem-0}
f^4(x,y,z)
&=\int_{\sigma}^z \frac{d}{d\xi} \left(f^4(x,y,\xi)\right) d\xi + f^4(x,y,\sigma) \notag\\
&=4\int_{\sigma}^z f^3(x,y,\xi)  f_{\xi}(x,y,\xi) d\xi + f^4(x,y,\sigma) \notag\\
&\leq 4\int_0^1 |f(x,y,\xi)|^3 |f_{\xi}(x,y,\xi)| d\xi +  f^4(x,y,\sigma).
\end{align}
Integrating (\ref{lem-0}) with respect to $\sigma$ over $[0,1]$, we obtain
\begin{align*}  
f^4(x,y,z) \leq    4\int_0^1 |f(x,y,\xi)|^3 |f_{\xi}(x,y,\xi)| d\xi + \int_0^1 f^4(x,y,\sigma) d\sigma,
\end{align*}
and by Cauchy-Schwarz inequality, we have
\begin{align}     \label{lem-1}
f^4(x,y,z)   \leq     4   \|f\|^3_{L^6_z}  \|f_z\|_{L^2_z}     + \|f\|_{L^4_z}^4.
\end{align}
Here we denote
$$\|f\|_{L^p_z}= \left(\int_0^1 |f(x,y,z)|^p dz\right)^{1/p}.$$
Now, if we denote $\|f\|_{L_z^{\infty}} =  \sup_{z\in[0,1]}  |f(x,y,z)| $, then the inequality (\ref{lem-1}) can be written as
\begin{align}   \label{lem-1'}
\|f\|_{L^{\infty}_z}   \leq C \|f\|^{3/4}_{L^6_z}  \|f_z\|^{1/4}_{L^2_z}     + \|f\|_{L^4_z}.
\end{align}

Thanks to the H\"older's inequality and (\ref{lem-1'}), we have
\begin{align}   \label{lem-2} 
&\int_{\Omega} |fgh| dx dy dz     \notag\\
&\leq   \int_{[0,L]^2} \|f\|_{L^{\infty}_z}   \|g\|_{L^2_z}   \|h\|_{L^2_z}  dx dy \notag\\
&\leq   C\int_{[0,L]^2} \left( \|f\|^{3/4}_{L^6_z}  \|f_z\|_{L^2_z}^{1/4}     + \|f\|_{L^4_z} \right)  \|g\|_{L^2_z}  \|h\|_{L^2_z} dx dy \notag\\
&\leq   C  \left[ \left( \int_{[0,L]^2}   \|f\|^3_{L^6_z}  \|f_z\|_{L^2_z} dx dy\right)^{1/4} + \|f\|_4 \right]   
\left(\int_{[0,L]^2} \|g\|_{L^2_z}^4 dx dy\right)^{1/4}  \|h\|_2    \notag\\
&\leq C\left(\|f\|_6^{3/4}  \|f_z\|_2^{1/4} +  \|f\|_4  \right) \left(\int_{[0,L]^2} \|g\|_{L^2_z}^4 dx dy\right)^{1/4}   \|h\|_2   .
\end{align}

Recall the Ladyzhenskaya inequality (see, e.g., \cite{Lady}) in the three-dimensional periodic domain $\Omega$:
\begin{align}   \label{iden-5}
\|\varphi\|_p \leq C_p \|\varphi\|_2^{\frac{6-p}{2p}} \left(\|\varphi\|_2 + \|\varphi_x\|_2\right)^{\frac{p-2}{2p}}  
\left(\|\varphi\|_2 + \|\varphi_y\|_2\right)^{\frac{p-2}{2p}} \left(\|\varphi\|_2 + \|\varphi_z\|_2\right)^{\frac{p-2}{2p}},
\end{align}
for $\varphi\in H^1(\Omega)$ and $p\in [2,6]$. By (\ref{iden-5}), one has
\begin{align}  \label{lem-3}
\|f\|_6 \leq C\left(\|f\|_2+\|\nabla_h f\|_2 \right)^{2/3} \left(\|f\|_2+\|f_z\|_2\right)^{1/3},
\end{align}
and 
\begin{align}  \label{lem-4}
\|f\|_4 \leq C\|f\|_2^{1/4}   \left(\|f\|_2+\|\nabla_h f\|_2 \right)^{1/2}   \left(\|f\|_2+\|f_z\|_2\right)^{1/4}.
\end{align}
By virtue of (\ref{lem-3}) and (\ref{lem-4}), we have
\begin{align}  \label{lem-5}
\|f\|_6^{3/4}  \|f_z\|_2^{1/4} +  \|f\|_4 \leq C  \left(\|f\|_2+\|\nabla_h f\|_2 \right)^{1/2} \left(\|f\|_2+\|f_z\|_2\right)^{1/2}.
\end{align}

Recall the Agmon's inequality (see, e.g., \cite{Agmon}) in one dimension: 
\begin{align}    \label{lem-6}
\|\phi\|_{L^{\infty}([0,L])} \leq C \|\phi\|_{L^2([0,L])}^{1/2} \|\phi\|_{H^1([0,L])}^{1/2}.
\end{align}
By using (\ref{lem-6}), we obtain 
\begin{align}   \label{lem-7}
&\int_{[0,L]^2} \|g\|_{L^2_z}^4 dx dy = \int_{[0,L]^2} \left(\int_0^1 g^2 dz\right)  \left(\int_0^1 g^2 dz\right)  dx dy \notag\\
& \leq C \left[\int_0^L \int_0^1 \left(\int_0^L g^2 dx\right)^{\frac{1}{2}}  \left(\int_0^L \left(g^2+g_x^2\right) dx\right)^{\frac{1}{2}} dz dy\right]  \notag\\
&\hspace{0.8 in}  \cdot \left[\int_0^L \int_0^1 \left(\int_0^L g^2 dy\right)^{\frac{1}{2}}  \left(\int_0^L \left(g^2+g_y^2\right) dy\right)^{\frac{1}{2}} dz dx\right]   \notag\\
&\leq C \|g\|_2^2 \left(\|g\|_2^2 + \|\nabla_h g\|_2^2\right).
\end{align}
By combining (\ref{lem-2}), (\ref{lem-5}) and (\ref{lem-7}), we conclude the proof of Lemma \ref{lemma1}.
\end{proof}

\vspace{0.3 in}

\noindent {\bf Acknowledgment.} The work of E.S.T. was supported in part by the ONR grant N00014-15-1-2333.

\vspace{0.3 in}

\bibliographystyle{amsplain}

\end{document}